\documentstyle[amscd,amssymb,verbatim]{amsart}
\newcommand{\Ob}{\operatorname{Ob}}

\newcommand{\ti}{\tilde}
\newcommand{\pr}{\operatorname{pr}}

\newcommand{\DD}{{\cal D}}

\newcommand{\lan}{\langle}
\newcommand{\ran}{\rangle}

\newcommand{\CC}{{\cal C}}

\newcommand{\si}{\sigma}

\newcommand{\ga}{\gamma}
\newcommand{\de}{\delta}

\newcommand{\D}{{\cal D}}
\numberwithin{equation}{subsection}

\newtheorem{thm}{Theorem}[section]

\newtheorem{lem}[thm]{Lemma}

\newcommand{\Pf}{\noindent {\it Proof}}
\newcommand{\id}{\operatorname{id}}

\newcommand{\ov}{\overline}
\newcommand{\we}{\wedge}
\renewcommand{\Im}{\operatorname{Im}}

\newcommand{\ra}{\rightarrow}
\renewcommand{\AA}{{\cal A}}

\newcommand{\HH}{{\cal H}}

\newcommand{\VV}{{\cal V}}

\renewcommand{\O}{{\cal O}}
\newcommand{\Om}{\Omega}
\newcommand{\dbar}{\overline{\partial}}
\newcommand{\Hom}{\operatorname{Hom}}
\newcommand{\Ext}{\operatorname{Ext}}

\renewcommand{\a}{\alpha}
\renewcommand{\b}{\beta}

\newcommand{\De}{\Delta}
\newcommand{\la}{\lambda}
\newcommand{\th}{\theta}
\newcommand{\C}{{\Bbb C}}
\newcommand{\R}{{\Bbb R}}
\newcommand{\Z}{{\Bbb Z}}
\newcommand{\Q}{{\Bbb Q}}
\newcommand{\La}{\Lambda}

\newcommand{\wt}{\widetilde}

\newcommand{\ed}{\qed\vspace{3mm}}

\title{Homological mirror symmetry with higher products}
\author{A. Polishchuk}

\begin{document}
\maketitle

\bigskip

The homological mirror symmetry conjecture formulated by M.~Kontsevich
in \cite{K}
claims that derived categories of Fukaya's symplectic $A_{\infty}$-categogy
$F(M)$ of a Calabi-Yau manifold $M$
and of coherent sheaves on a mirror dual Calabi-Yau manifold $X$
are equivalent.
In particular, this means that one can identify the associative product on
$\Ext$-groups between coherent sheaves on $X$ with the corresponding product
in the Floer cohomology of Lagrangians submanifolds in $M$ (defined by
Fukaya in \cite{F1}).
The drawback of this conjecture is that one has an $A_{\infty}$-category
on the symplectic side of the story and the usual category on
the complex side, so one has to make the usual category out of $F(M)$.
In this note we fix this problem by constructing
an $A_{\infty}$-category on the complex side and formulate a more
general conjecture involving $A_{\infty}$-categories on both sides.

Let $X$ be a compact complex manifold equipped with a hermitian metric.
Inspired by Merkulov's paper \cite{M}
we define an $A_{\infty}$-category $\DD^b_{\infty}(X)$ which
is a refined version of the derived category of $X$.
The objects of $\DD^b_{\infty}(X)$ are bounded complexes of
holomorphic vector bundles on $X$ equipped with hermitian metrics.
The morphisms from $E^{\bullet}$ to $F^{\bullet}$ are elements of
$\Ext(E^{\bullet},F^{\bullet})=\oplus_q
H^q(X,(E^{\bullet})^*\otimes F^{\bullet})$ which can be thought of as
harmonic $(0,q)$-forms with values in $(E^{\bullet})^*\otimes F^{\bullet}$.
The $A_{\infty}$-structure
has $m_1=0$ while $m_2$ is the usual composition of $\Ext$'s. In particular,
forgetting higher products we obtain the usual derived category of $X$.
The higher products measure in some sense to which extent the product
of harmonic forms with values in hermitian bundles fails to be
harmonic. More precisely, we use the construction of $A_{\infty}$-structure
on a subcomplex homotopy equivalent to a dg-algebra (see \cite{GS}, 
\cite{GLS}, \cite{M}). This construction gives an $A_{\infty}$-structure
on the algebra $\Ext(E,E)$ where $E$ is a hermitian holomorphic
vector bundle (or complex of such bundles).
It follows from the results of the homological perturbation theory 
(see \cite{GLS}) that up to homotopy equivalence
our $A_{\infty}$-structures do not depend on choices of metrics.

Now the more natural formulation the homological mirror
conjecture is that for mirror dual Calabi-Yau manifolds 
there exists an $A_{\infty}$-functor from $F(M)$ to
$\DD^b_{\infty}(X)$ which is a homotopy equivalence on morphisms.
M.~Kontsevich informed the author that such a formulation is not new
provided that one replaces $\DD^b_{\infty}(X)$ by a dg-version of
the derived category constructed either using Cech resolutions or
using Dolbeault resolutions (see e.g. \cite{BK}). The dg-category obtained
from Dolbeault resolutions is $A_{\infty}$-equivalent to $\DD^b_{\infty}(X)$.
The only advantage of $\D^b_{\infty}(X)$ is that it has 
finite-dimensional $\Hom$'s and the higher products satisfy an explicit
cyclic symmetry as we show in \ref{Serresec}.

So far the only non-trivial case in which the original
Kontsevich's conjecture is  verified
is that of elliptic curves considered in \cite{PZ}.
In the present paper we show that generic
triple products on the subcategory of line bundles on elliptic curves
defined using complex structure and using symplectic structure are
canonically homotopic to each other. This gives a partial verification
for our generalized version of the homological mirror conjecture in the case
of elliptic curves.

\section{$A_{\infty}$-categories}

\subsection{}
Recall that a ($\Z$-graded) $A_{\infty}$-algebra is a $\Z$-graded
vector space $A$
equipped with linear maps $m_k: A^{\otimes k}\ra A$ for
$k\ge 1$ of degree $2-k$ satisfying the following $A_{\infty}$-constraints
\begin{align*}
&\sum_{k+l=n+1}\sum_{j=0}^{k-1}
(-1)^{l(\ti{a}_1+\ldots+\ti{a}_j)+j(l-1)+(k-1)l}
m_k(a_1,\ldots,a_j,m_l(a_{j+1},\ldots,a_{j+l}),a_{j+l+1},\ldots,a_n)\\
&=0
\end{align*}
where $\ti{a}_i=\deg(a_i)\mod(2)$.
A more conceptual way to spell this axiom is to say that $m_k$
are components of a coderivation $\partial$ on the bar-construction of $A$
such that $\partial^2=0$ (see \cite{JS}).
For a pair of $A_{\infty}$-algebras $A$ and $B$ there is a natural
notion of a $A_{\infty}$-morphism from $A$ to $B$. Namely,
such a morphism consists of the data $(f_n, n\ge 1)$ where
$f_n:A^{\otimes n}\ra B$ is a linear of degree $1-n$ such that
\begin{align*}
&\sum_{1\le k_1<k_2<\ldots<k_i=n}\pm
m_i(f_{k_1}(a_1,\ldots,a_{k_1}),f_{k_2-k_1}(a_{k_1+1},\ldots,a_{k_2}),
\ldots,f_{n-k_{i-1}}(a_{k_{i-1}+1},\ldots,a_n))\\
&=\sum_{k+l=n+1}\sum_{j=0}^{k-1}\pm
f_k(a_1,\ldots,a_j,m_l(a_{j+1},\ldots,a_{j+l}),a_{j+l+1},\ldots,a_n)
\end{align*}
Again this is equivalent to having a morphism between the corresponding
bar-constructions in the category of differential coalgebras.

One can compose $A_{\infty}$-morphisms in the natural way.
The identity $A_{\infty}$-morphism consists of $f_1=\id$,
$f_n=0$ for $n\ge 2$.
If $(f_n):A\ra B$ is an $A_{\infty}$-morphism such
that $f_1$ is an isomorphism of underlying abelian groups then
$f_1^{-1}$ extends to an $A_{\infty}$-morphism $B\ra A$ which
is inverse to $(f_n)$. In the case $B$ and $A$ have the same
underlying spaces and $f_1=\id$ we will call
sometimes the data $(f_n,n\ge 2)$ a homotopy between two
$A_{\infty}$-structures $(m_n)$ and $(m'_n)$ on the same space.
For example, if we have $m_1=m'_1=0$ and $m_2(a,b)=m'_2(a,b)=ab$
then $f_2$ gives a homotopy between $m_3$ and $m'_3$ in the following sense:
$$m_3'(a_1,a_2,a_3)-m_3(a_1,a_2,a_3)=
(-1)^{\wt{a}_1}a_1f_2(a_2,a_3)-f_2(a_1,a_2)a_3-f_2(a_1a_2,a_3)+
f_2(a_1,a_2a_3).$$

\subsection{}
The definition of $A_{\infty}$-category is very similar to that
of an $A_{\infty}$-algebras (see \cite{F1}). 
Namely, an $A_{\infty}$-category $\CC$
consists of a class of objects $\Ob\CC$,
for every pair of objects $E_1$ and $E_2$ graded spaces of morphisms
$\Hom(E_1,E_2)$ equipped with a differential $m_1$, and compositions
$$m_k:\Hom(E_1,E_2)\otimes\Hom(E_2,E_3)\otimes\ldots\otimes
\Hom^*(E_k,E_{k+1})\ra\Hom(E_1,E_{k+1})$$
of degree $2-k$ for all $k\ge2$. The associativity constraint is
that these compositions define a structure of $A_{\infty}$-algebra
on $\oplus_{ij}\Hom(E_i,E_j)$ for every collection
$E_1,\ldots,E_n\in\Ob\CC$.

An $A_{\infty}$-functor (see \cite{F2}) $\phi:\CC\ra\CC'$ between
$A_{\infty}$-categories consists of the map $\phi:\Ob\CC\ra\Ob\CC'$ and
of the linear maps
$$f_k:\Hom_{\CC}(E_1,E_2)\otimes\Hom_{\CC}(E_2,E_3)\otimes\ldots\otimes
\Hom_{\CC}(E_k,E_{k+1})\ra\Hom_{\CC'}(\phi(E_1),\phi(E_{k+1}))$$
of degree $1-k$ for $k\ge 1$,
which define $A_{\infty}$-morphisms
$\oplus_{ij}\Hom_{\CC}(E_i,E_j)\ra\oplus_{ij}\Hom_{\CC'}(\phi(E_i),
\phi(E_j))$.


\subsection{}
Recall that if $(A,d)$ is a dg-algebra, $Q:A\ra A$ is an odd operator
then the subcomplex $B=(1-dQ-Qd)(A)\subset A$ has a canonical
$A_{\infty}$-structure defined as follows (see \cite{M}).
The differential $m_1$ is just the restriction of $d$ to $B$.
Then one defines inductively linear maps $\la_n:A^{\otimes n}\ra A$
for $n\ge 2$ by setting $\la_2(a_1,a_2)=a_1a_2$,
\begin{align*}
&\la_n(a_1,\ldots,a_n)=(-1)^{n-1}
(Q\la_{n-1}(a_1,\ldots,a_{n-1})a_n-(-1)^{n\ti{a}_1}
a_1(Q\la_{n-1}(a_2,\ldots,a_n))-\\
&\sum_{k+l=n, \atop k,l\ge 2}(-1)^{k+(l-1)(\ti{a}_1+\ldots+\ti{a}_k)}
(Q\la_k(a_1,\ldots,a_k))(Q\la_l(a_{k+1},\ldots,a_n)).
\end{align*}
Now the products $m_n$ for $n\ge 2$  are defined by the formula
$$m_n(b_1,\ldots,b_n)=\pr(\la_n(b_1,\ldots,b_n))$$
where $\pr=1-Qd-dQ$.
For example, $m_2(b_1,b_2)=\pr(b_1b_2)$,
$$m_3(b_1,b_2,b_3)=\pr(Q(b_1b_2)b_3-(-1)^{\ti{b}_1}b_1Q(b_2b_3)).$$

Assume in addition that we have $\pr|_B=\id$, $\pr\circ Q=Q^2=0$.
Then we can apply the homological perturbation theory
(see \cite{GLS}, 4.2) to conclude that the natural embedding
$B\ra A$ and the projection $\pr:A\ra B$ extend to $A_{\infty}$-morphisms
of $A_{\infty}$-algebras (where $A$ has trivial $m_k$ with $k\ge 3$).
In the situation below we will also have $d|_B=0$, so that $B$ is
isomorphic to the cohomology $H(A)$ of $(A,d)$. Thus, in this
particular case the perturbation theory implies that $A_{\infty}$-structures
on $H(A)$ obtained by different operators $Q$ are homotopic.
Note that the higher products on $H(A)$ thus defined coincide with
the Massey products when the latter are well-defined and univalued.

We apply this machinery in the following situation:
$(A,d)=(C^{\infty}(X,\Om^{0,\bullet}(E^*\otimes E),\dbar))$
is the Dolbeault complex computing $\Ext(E,E)$, where
$E$ is a holomorphic vector bundle on a compact complex manifold $X$.
Now assume that $X$ and $E$ are equipped with hermitian metrics. Then
there is an induced metric on $E^*\otimes E$ which defines a
conjugate-linear $*$-operator
$$\ov{*}:\Om^{p,q}(E^*\otimes E)\ra\Om^{n-p,n-q}(E^*\otimes E).$$
The induced metric on the Dolbealt complex is given by
$$(\a,\b)=\int_X\a\we\ov{*}\b.$$
The conjugate operator to $\dbar$ with respect to this metric is
$$\dbar^*=-\ov{*}\dbar\ov{*}$$
(see \cite{W}).
Now we set $Q=\dbar^*G$ where
$G$ is the Green operator corresponding to the
$\dbar$-laplacian $\De=\dbar^*\dbar+\dbar\dbar^*$. Then
by definition 
$$\pr=1-\dbar Q-Q\dbar=1-\De G-G\De$$
is the projector on the space $\HH^{0,\bullet}(E^*\otimes E)$
of harmonic forms of type $(0,\bullet)$.
Thus, the above general construction gives an $A_{\infty}$-structure
on $\HH^{0,\bullet}(E^*\otimes E)\simeq\Ext(E,E)$
extending the standard algebra structure on it (and having $m_1=0$).

Similarly, fixing a hermitian metric on $X$
we can define an $A_{\infty}$-category $\VV_{\infty}(X)$
as follows.
The objects of $\VV_{\infty}(X)$ are hermitian vector bundles on $X$.
Morhisms from $E$ to $F$ in $\CC(X)$ are elements of
$\Ext(E,F)$ (identified with
the space of harmonic forms $\HH^{0,\bullet}(E^*\otimes F)$).
The composition
\begin{equation}\label{mk}
m_k:\Ext(E_0,E_1)\otimes\ldots\otimes\Ext(E_{k-1},E_k)\ra
\Ext(E_0,E_k)
\end{equation}
is defined as follows.
Consider the hermitian holomorphic vector bundle
$E=\oplus_{i=0}^k E_i$. Then the left hand side of (\ref{mk})
is the direct summand of $\Ext(E,E)^{\otimes k}$ while
the right hand side is the direct summand of
$\Ext(E,E)$. So we just define (\ref{mk}) as the corresponding
component of the map
$$m_k:\Ext(E,E)^{\otimes k}\ra\Ext(E,E)$$
coming from the $A_{\infty}$-structure on $\Ext(E,E)$.


\subsection{} The definition of the category $\DD^b_{\infty}(X)$
for a hermitian complex manifold $X$ follows
the same pattern. One just has to observe that there is an analogue
of the theory of harmonic forms for a complex of
hermitian holomorphic vector bundles $(E^{\bullet},\de)$ (the differential
$\de$ is holomorphic). Namely, we consider
the complex of forms $\Om^{0,\bullet}(E^{\bullet})=\oplus_{p,q}
\Om^{0,p}(E^q)$ with the total differential $\dbar+\de$ and total grading
$p+q$. Then the operator
$\dbar+\de$ is elliptic, so we can introduce the corresponding Laplace
operator $\Delta_{\dbar+\de}$, Green operator $G_{\dbar+\de}$, etc.
In particular, we have a notion of harmonic forms with values in
$E^{\bullet}$. They don't have to be bihomogeneous so the space
of harmonic forms has just one grading coming from the total grading
on $\Om^{0,\bullet}(E^{\bullet})$. Now we can apply the same procedure as
before to define the $A_{\infty}$-structure on the space of harmonic
forms which can be identified with
$\Ext(E^{\bullet},E^{\bullet})$ (=morphisms in the derived category
$\D^b(X)$). Now $\DD^b_{\infty}(X)$
has as objects complexes of hermitian holomorphic vector bundles with
morphisms between $E^{\bullet}$ and $F^{\bullet}$ being the space
of $\Ext(E^{\bullet},F^{\bullet})$.
The $A_{\infty}$-structure is defined as above.

\subsection{} One can generalize slightly the above construction
as follows. For a holomorphic vector bundle
$E$ on $X$ let us denote $H^{p,q}(E)=H^q(X,\Om^p_{hol}\otimes E)$. 
Then for an $n$-tuple $E_1,\ldots,E_n$ of hermitian bundles on $X$ 
(a hermitian metric on $X$ is fixed) we can define products
$$m_n:H^{p_1,q_1}(E_1)\otimes\ldots\otimes H^{p_n,q_n}(E_n)\ra
H^{p_1+\ldots+p_n,q_1+\ldots+q_n-n+2}(E_1\otimes\ldots\otimes E_n).$$
Indeed, this product is obtained by applying the general construction
to the Dolbeault complex of forms of all types with values in the  
algebra $\AA=\oplus_{i\le j}\AA_{ij}$, where
$\AA_{ij}=E_i\otimes E_{i+1}\ldots\otimes E_j$,
equipped with the operator $Q=\dbar^*G$.
The products (\ref{mk}) can be recovered from this more general construction.
The further possible generalization is to let the product $m_k$
depend on a collection of hermitian metrics $h_{ij}$, $1\le i<j\le k$ on $X$,
by considering the operator $Q$ on the Dolbeault complex with values in
$\AA$ such that $Q=\oplus Q_{ij}$, where $Q_{ij}$ is defined using the
metric $h_{ij}$ (note that only $Q_{ij}$ with $i<j$ are used in the definition
of $m_k$). These generalized products still satisfy an analogue of
the $A_{\infty}$-constraint.

\subsection{}\label{Serresec} 
The higher products defined above are compatible
with the Serre duality in the way described below. Recall 
that the Serre duality is
a perfect pairing between $H^{p,q}(E)$ and $H^{d-p,d-q}(E^*)$.
where $d$ is the complex dimension of $X$.
If $X$ and $E$ are equipped with hermitian metrics then this pairing is
induced by the pairing
$$\lan \a,\b\ran=\int \a\we\b$$
between the corresponding Dolbeault complexes. More generally,
for any bilinear holomorphic pairing $b:E\otimes_{\O} F\ra\O$ we have the
induced pairing $\lan\cdot,\cdot\ran_b$ between
$H^{p,q}(E)$ and $H^{d-p,d-q}(F)$.
Now assume that we have $n+1$ hermitian holomorphic bundles
$E_1,\ldots, E_{n+1}$ and a polylinear holomorphic map
$$f:E_1\otimes_{\O}\ldots\otimes_{\O} E_{n+1}\ra\O.$$

\begin{thm}\label{Serre}
One has
$$\lan m_n(\a_1,\ldots,\a_n),\a_{n+1}\ran_f=(-1)^{n(\ti{\a}_1+1)}
\lan \a_1, m_n(\a_2,\ldots,\a_{n+1})\ran_f$$
where $\a_i\in H^{p_i,q_i}(E_i)$, $i=1,\ldots,n+1$,
$p_1+\ldots+p_{n+1}=d$, $q_1+\ldots+q_{n+1}=d+n-2$.
\end{thm}

\begin{lem}\label{Slem}
Let $E$ be a hermitian holomorphic bundle. Then for any
$\a\in\Om^{p,q+1}(E)$, $\b\in\Om^{d-p,d-q}(E^*)$ one has
$$\int_X Q\a\we\b=(-1)^{\ti{\a}}\int_X \a\we Q\b$$
where $Q=\dbar^* G$.
\end{lem}

\Pf . This follows from self-adjointness of $G$ and from the formula
$$\int\dbar^*\a\we\b=(-1)^{p+q+1}\int\a\we\dbar^*\b$$
which in turn follows from the equalities $\dbar^*=-\ov{*}\dbar\ov{*}$
and $\b=(-1)^{p+q}\ov{*}(\ov{*}\b)$.
\ed

\noindent
{\it Proof of theorem \ref{Serre}}.
First of all Lemma \ref{Slem} implies that
for a harmonic form $\b$ one has
$$\int\pr(\a)\we\b=\int\a\we\b.$$
Thus, we have to prove the identity
\begin{equation}\label{Seq}
\int\la_n(\a_1,\ldots,\a_n)\we\a_{n+1}=(-1)^{n(\ti{\a}_1+1)}
\int\a_1\we\la_n(\a_2,\ldots,\a_{n+1}).
\end{equation}
Following \cite{M} let us rewrite the inductive definition of $\la_n$
as follows:
$$\la_{n}(x_1,\ldots,x_n)=\sum_{k+l=n,\atop k,l\ge1}(-1)^{k+(l-1)
(\ti{x}_1+\ldots+\ti{x}_k)}
Q\la_k(x_1,\ldots,x_k) Q\la_l(x_{k+1},\ldots,n)$$
where we denote formally $\la_1=-Q^{-1}$ so that
$Q\la_1(x)=-x$.
Hence, the LHS of (\ref{Seq}) is equal to the integral of
\begin{align*}
&\sum_{k+l=n,\atop k,l\ge1}(-1)^{k+(l-1)(\ti{x}_1+\ldots+\ti{x}_k)}
Q\la_k(\a_1,\ldots,\a_k) Q\la_l(\a_{k+1},\ldots,\a_n)\a_{n+1}=\\
&(-1)^{1+n\ti{\a}_1}\a_1 Q\la_{n-1}(\a_2,\ldots,\a_n)\a_{n+1}+\\
&(-1)^{n}\sum_{k+l=n+1,\atop k,l\ge 2}(-1)^{l(\ti{\a}_1+\ldots+\ti{\a}_k)}
Q\la_k(\a_1,\ldots,\a_k) \la_l(\a_{k+1},\ldots,\a_{n+1})-\\
&(-1)^{n}\sum_{k+l+m=n+1,\atop k,m\ge2;l\ge1}(-1)^{l+(m-1)(\ti{\a}_{k+1}+
\ldots+\ti{\a}_{k+l})+(l+m)(\ti{\a}_1+\ldots+\ti{\a}_k)}
Q\la_k(\a_1,\ldots,\a_k)\times\\
&Q\la_l(\a_{k+1},\ldots,\a_{k+l})Q\la_m(\a_{k+l+1},\ldots,\a_{n+1}).
\end{align*}
Similarly, we can rewrite the RHS of (\ref{Seq}) as the integral of
\begin{align*}
&(-1)^{1+n(\ti{\a}_1)}\a_1 Q\la_{n-1}(\a_2,\ldots,\a_n)\a_{n+1}+\\
&(-1)^{n}\sum_{k+l=n+1,\atop k,l\ge 2}
(-1)^{k+(l-1)(\ti{\a}_1+\ldots+\ti{\a}_k)}
\la_k(\a_1,\ldots,\a_k) Q\la_l(\a_{k+1},\ldots,\a_{n+1})-\\
&(-1)^{n}\sum_{k+l+m=n+1,\atop k,m\ge2;l\ge1}(-1)^{l+(m-1)(\ti{\a}_{k+1}+
\ldots+\ti{\a}_{k+l})+(l+m)(\ti{\a}_1+\ldots+\ti{\a}_k)}
Q\la_k(\a_1,\ldots,\a_k)\times\\
&Q\la_l(\a_{k+1},\ldots,\a_{k+l})Q\la_m(\a_{k+l+1},\ldots,\a_{n+1}).
\end{align*}
It remains to apply Lemma \ref{Slem} to conclude that
\begin{align*}
&\int Q\la_k(\a_1,\ldots,\a_k) \la_l(\a_{k+1},\ldots,\a_{n+1})=\\
&(-1)^{k+\ti{\a}_1+\ldots+\ti{\a}_k}
\int\la_k(\a_1,\ldots,\a_k) Q\la_l(\a_{k+1},\ldots,\a_{n+1}).
\end{align*}
\ed

\section{Case of elliptic curve}\label{ell}
In this section we'll study the category $\VV_{\infty}(X)$ in the case
when $X=E_{\tau}=\C/\Z+\Z\tau$ is a complex elliptic curve.
All holomorphic vector bundles on elliptic curves are obtained
from line bundles using the operations of push-forward under isogeny,
tensoring with local systems, and direct sums. Therefore, it is
plausible that the entire $A_{\infty}$-category $\VV_{\infty}(E_{\tau})$
can be
reconstructed from its full $A_{\infty}$-subcategory consisting of
line bundles (see the argument of \cite{PZ} dealing with
the case of $m_2$). More precisely, we choose canonical 
representatives in each isomorphism class of holomorphic line bundles
on $E_{\tau}$ and canonical hermitian metrics on them and
consider the corresponding subcategory in $\VV_{\infty}(E_{\tau})$.
We show that generic triple products in this subcategory are homotopic
to the triple products between the corresponding objects
in the Fukaya category of the symplectic torus.
More precisely we consider {\it transversal} products, i.e.  
such that the corresponding configuration of geodesic circles
in $\R^2/\Z^2$ is transversal.

\subsection{}
Let $p:\C\ra E_{\tau}$ be the projection. We
denote by $L$ the holomorphic line bundle of degree 1 on $E_{\tau}$
equipped with a trivialization of $p^*L$ such that the 
theta-function
$$\theta(x,\tau)=\sum_{n\in\Z}\exp(\pi i\tau n^2+2\pi i n x)$$
is a pull-back of a holomorphic section of $L$.
For every $u\in\C$ let us denote by $L(k,u)$ the holomorphic line bundle
$L^{k-1}\otimes t^*_u L$ on $E_{\tau}$.
Up to an isomorphism $L(k,u)$ depends only on $u$ modulo the
lattice $\Z+\Z\tau$. More precisely, we have $L(k,u)=L(k,u+1)$
and a natural isomorphism
$$L(k,u)\wt{\ra} L(k,u+\tau):f(x)\mapsto \exp(-2\pi x)f(x).$$
We fix a hermitian
metric on $L(k,u)$ which under the natural trivialization of $\pi^*L(k,u)$
has form
$$\lan f, g\ran=\int_{\C/\Z+\Z\tau}f(x)\ov{g(x)}
\exp(-2\pi t(kx_2^2+2x_2u_2))dx_1dx_2$$
where $t=\Im(\tau)$, $x=x_1+\tau x_2$, $u=u_1+\tau u_2$.
The natural isomorphisms $L(k,u)\simeq L(k,u+m+n\tau)$ are compatible
with metrics. Also we have an isomorphism of hermitian line bundles
$$L(k,u)\otimes L(k',u')\simeq L(k+k',u+u').$$

Note that every $C^{\infty}$-section $f(x)$ of $L(k,u)$ can be written
in the form
$$f(x)=\sum_{n\in\Z}\varphi_n(x_2)\exp(2\pi i nx_1).$$
If $k\neq 0$ then $|k|$ functions $(\varphi_n, n=0,\ldots,|k|)$ determine the
rest. Moreover, the functions
$\varphi_n(x_2)\exp(\pi i(k\tau x_2^2+2ux_2))$
belong to the Schwarz space of functions on $\R$
with derivatives of all orders rapidly decreasing at infinity .
If $g(x)=\sum_{n\in\Z}\psi_n(x_2)\exp(2\pi i nx_1)$
is another $C^{\infty}$-section of $L(k,u)$ then we have
\begin{equation}\label{herm}
\lan f, g\ran=\sum_{n\in\Z/k\Z}\int_{\R}\varphi_n(s)\ov{\psi_n(s)}
\exp(-2\pi t(ks^2+2u_2s))ds.
\end{equation}

On the other hand, the line bundle $L(0,u)$ is
trivial as $C^{\infty}$-bundle.
Indeed, $C^{\infty}$-sections of $L(0,u)$ are functions $f(x)$ on $\C$
such that $f(x+1)=f(x)$, $f(x+\tau)=\exp(-2\pi i u)f(x)$. In other words,
the function $\exp(2\pi i ux_2)f(x)$, where $x=x_1+\tau x_2$,
is doubly periodic with respect
to $\Z+\Z\tau$. In particular, any $C^{\infty}$-section of $L(0,u)$
can be written as a Fourier series
$$f(x)=\sum_{m,n}c_{m,n}\exp(2\pi i(mx_1+ (n-u)x_2))$$
where $\varphi_{u,m,n}(x)=\exp(2\pi i(mx_1+(n-u)x_2))$ form an orthonormal
basis with respect to the hermitian metric on $L(0,u)$.
Note that $\dbar \varphi_{u,m,n}(x)=\frac{\pi}{t}(m\tau-n+u) \varphi_{u,m,n}(x)$.
Hence, if $c_{m,n}$ are the Fourier coefficients of
$f\in C^{\infty}(L(0,u))$
then $\frac{\pi}{t}(m\tau-n+u)c_{m,n}$ are the Fourier coefficients of
$\dbar f$.

\subsection{}
If $k>0$ then the basis of global sections of $L(k,u)$ is given
by $k$ theta-functions with characteristics
$$\theta_{a/k}(kx+u,k\tau); a\in\Z/k\Z$$
where for $r\in\Q/\Z$ we denote
$$\theta_r(x,\tau)=\sum_{n\in\Z}\exp(\pi i\tau(n+r)^2+2\pi i (n+r)x)
=\exp(\pi i \tau r^2+2\pi i r x)\th(x+r\tau,\tau).$$
Using (\ref{herm}) it is easy to see that
this basis is orthogonal with respect to the hermitian metric and 
$$||\theta_{a/k}(kx+u,k\tau)||^2=\frac{1}{\sqrt{2tk}}
\exp(\frac{2\pi t u_2^2}{k}).$$

Similarly, for $k>0$ the basis of harmonic $(0,1)$-forms with values
in $L(-k,-u)$ is given by
$$\ov{\theta_{a/k}(kx+u,k\tau)}\exp(-2\pi t(kx_2^2+2x_2u_2))d\ov{x},
a\in\Z/k\Z$$

Note that $\th_{b/k}(kx+u,k\tau)\ov{\th_{a/k}(kx+v,k\tau)}
\exp(-2\pi t(kx_2^2+2x_2v_2))$ is a section of $L(0,u-v)$.
The corresponding Fourier coefficients are
\begin{align*}
&c_{m,n}=\int_{\C/\Z+\Z\tau}
\th_{b/k}(kx+u,k\tau)\ov{\th_{a/k}(kx+v,k\tau)}\times\\
&\exp(-2\pi t(kx_2^2+2x_2v_2)-2\pi i(mx_1+(n-u+v)x_2))dx_1dx_2=\\
&\lan \th_{b/k}(kx+u,k\tau)\exp(-2\pi i(mx_1+(n-u+v)x_2)),
\th_{a/k}(kx+v,k\tau)\ran,
\end{align*}
where we take the hermitian product of sections of $L(k,v)$.
Using the formula (\ref{herm}) we find that
$c_{m,n}=0$ unless $m\equiv b-a\mod(k)$ in which case
\begin{align*}
&c_{m,n}=\exp(\frac{\pi i}{k}((m+a)^2\tau-a^2\ov{\tau}
+2(m+a)u-2a\ov{v}))\times\\
&\int_{\R}\exp(-2\pi t ks^2+
2\pi i s((m+a)\tau-a\ov{\tau}-n+u-\ov{v}))ds=\\
&\frac{1}{\sqrt{2t k}}\exp(-\frac{\pi}{2t k}
(|\ga|^2+2\ov{\ga}u-2\ga\ov{v}+(u-\ov{v})^2)
+\frac{\pi i n}{k}(m+2a)). 
\end{align*}
where $\ga=m\tau-n$.

\subsection{}\label{comp} We start by computing the triple product
$$m_3:H^1(L(-k,0))\otimes H^0(L(k,u))\otimes H^0(L(l,v)\ra
H^0(L(l,u+v))$$
where $k,l>0$, $u\not\in\Z+\Z\tau$.
By definition this product is given by the formula
$$m_3(\a,\b_1,\b_2)=\pr(Q(\a\b_1)\b_2+\a Q(\b_1\b_2))$$
where $\a,\b_1,\b_2$ are harmonic representatives of the corresponding classes,
$\pr$ is the harmonic projector.
Note that $\b_1\b_2$ is a holomorphic section of $L(k+l,u+v)$, hence
$Q(\b_1\b_2)=0$. On the other hand, since $u\not\in\Z+\Z\tau$ we have
$H^*(L(0,u))=0$, hence, $Q(\a\b_1)=\dbar^{-1}(\a\b_1)$, so
$m_3(\a,\b_1,\b_2)=\pr(\dbar^{-1}(\a\b_1)\b_2)$.
We can compute this using explicit bases of theta functions: 
\begin{align*}
&m_3(\ov{\th_{a/k}(kx,k\tau)}\exp(-2\pi t kx_2^2)d\ov{x},
\th_{b/k}(kx+u,k\tau),\th_{c/l}(lx+v,l\tau))=
\sqrt{2tl}\times \\
&\sum_{d=0}^{l-1}
\lan F(x)\th_{c/l}(lx+v,l\tau),\th_{d/l}(lx+u+v,l\tau)\ran
\exp(-\frac{2\pi t(u_2+v_2)^2}{l})\th_{d/l}(lx+u+v,l\tau)
\end{align*}
where $F(x)\in C^{\infty}(L(0,u))$ is determined by
$$\dbar F(x)=
\ov{\th_{a/k}(kx,k\tau)}\exp(-2\pi tkx_2^2)\th_{b/k}(kx+u,k\tau).$$
Using the computation of the Fourier coefficients of the right hand side
we find that
$$F(x)=\sum_{m,n,\atop m\equiv b-a\mod(k)}a_{m,n}\varphi_{u,m,n}(x)$$
where
$$a_{m,n}=
\frac{t}{\pi\sqrt{2tk}(m\tau-n+u)}\exp(-\frac{\pi}{2tk}
(|m\tau-n|^2+2(m\ov{\tau}-n)u+u^2)
+\frac{\pi i}{k}(mn+2na)). $$

We have
\begin{align*}
&\lan F(x)\th_{c/l}(lx+v,l\tau),\th_{d/l}(lx+u+v,l\tau)\ran=\\
&\int_{\C/\Z+\Z\tau}
F(x)\th_{c/l}(lx+v,l\tau)\ov{\th_{d/l}(lx+u+v,l\tau)}
\exp(-2\pi t(lx_2^2+2x_2(u_2+v_2)))dx_1dx_2=\\
&\sum_{m,n,\atop m\equiv b-a\mod(k)}a_{m,n}\ov{b_{m,n}}
\end{align*}
where $b_{m,n}$ are the Fourier coefficients of
$\th_{d/l}(lx+u+v,l\tau)\ov{\th_{c/l}(lx+v,l\tau)}
\exp(-2\pi t(lx_2^2+2x_2v_2))$ (this is a $C^{\infty}$-section
of $L(0,u)$). We have $b_{m,n}=0$ unless $m\equiv d-c\mod(l)$
in which case
$$b_{m,n}=\frac{1}{\sqrt{2tl}}\exp(-\frac{\pi}{2tl}
(|m\tau-n|^2+2(m\ov{\tau}-n)(u+v)-2(m\tau-n)\ov{v}+(u+v-\ov{v})^2)
+\frac{\pi i n}{l}(m+2c)).$$
Denoting $w=u+v$ we obtain
\begin{align*}
&\lan F(x)\th_{c/l}(lx+v,l\tau),\th_{d/l}(lx+u+v,l\tau)\ran=
\frac{1}{2\pi\sqrt{kl}}\times\\
&\sum
\frac{\exp(-\frac{\pi(k+l)}{2tkl}Q(\ga,u)
+\frac{2\pi i}{l}((u+n)w_2-mw_1)+\frac{2\pi t}{l}w_2^2+
\pi i n(\frac{m+2c}{l}-\frac{m+2a}{k}))}{\ga+u}
\end{align*}
where the sum is taken over 
$\ga=m\tau+n\in\Z+\Z\tau$ such that $m\equiv b-a\mod(k)$ and 
$m\equiv d-c\mod(l)$, and
$$
Q(\ga,u)=|\ga|^2+2\ov{\ga}u+u^2
$$
Thus, we get the following formula for the triple product:
$$m_3(\a,\b_1,\b_2)=\sum_{d\in\Z/d\Z}G_{a,b,c}^d(u,w)
\th_{d/l}(lx+w,l\tau)$$
where
$\a=\frac{\pi \sqrt{2 k}}{\sqrt{t}}\ov{\th_{a/k}(kx,k\tau)}
\exp(-2\pi t kx_2^2)d\ov{x}\in H^1(L(-k,0))$,
$\b_1=\th_{b/k}(kx+u,k\tau)\in H^0(L(k,u))$,
$\b_2=\th_{c/l}(lx+v,l\tau))\in H^0(L(l,v))$,
\begin{align*}
&G_{a,b,c}^d(u,w)=\\
&\sum_{\ga=m\tau+n,\atop m\equiv b-a(k), m\equiv d-c(l)}
\frac{\exp(-\frac{\pi(k+l)}{2t kl}Q(\ga,u)
+\frac{2\pi i}{l}((u+n)w_2-mw_1)+
\pi i n(\frac{m+2c}{l}-\frac{m+2a}{k})}{\ga+u}
\end{align*}
The  function $G_{a,b,c}^d(u,w)$ is meromorphic in $u$ with poles of order 1
at the lattice points and satisfies
\begin{equation}\label{per1}
G_{a,b,c}^d(u+1,w)=\exp(2\pi i(\frac{b}{k}-\frac{c}{l}))G_{a,b,c}^d(u,w)
\end{equation}
\begin{equation}\label{per2}
G_{a,b,c}^d(u+\tau,w)=\exp(-\frac{\pi i(k+l)}{kl}(2u+\tau)+2\pi i\frac{w}{l})
G_{a,b+1,c-1}^d(u,w)
\end{equation}

\subsection{} Now let us compare the above triple product with
the corresponding product in the Fukaya category. 
According to the dictionary from \cite{PZ} the corresponding four
special lagrangians in the symplectic torus $\R^2/\Z^2$ are
$\La_1=(t,0)$, $\La_2=(t,-kt)$, $\La_3=(t,-u_2)$,
$\La_4=(t,lt-w_2)$. The lagrangians $\La_1$ and $\La_2$ are
equipped with trivial connections, while $\La_3$ (resp. $\La_4$) comes with
the connection $-2\pi iu_1dx$ (resp. $-2\pi iw_1dx$). The morphisms
$\a$, $\b_1$ and $\b_2$ in the derived category of $E_{\tau}$ correspond
to the intersection points
$(\frac{a}{k},0)\in\La_1\cap\La_2$, $(\frac{u_2+b}{k},-u_2)\in\La_2\cap\La_3$,
and $(\frac{w_2-u_2+c}{l},-u_2)\in\La_3\cap\La_4$ respectively, while
$\th_{d/l}(lx+w,l\tau)$ corresponds to the point
$(\frac{w_2+d}{l},0)\in\La_4\cap\La_1$.
Up to some factors which appear in the definition of the functor in \cite{PZ}
the corresponding Fukaya product is given by
$$m'_3(\a,\b_1,\b_2)=\sum_{d\in\Z/d\Z}F_{a,b,c}^d(u,w)
\th_{d/l}(lx+w,l\tau)$$
where
\begin{align*}
&F_{a,b,c}^d(u,w)=\\
&-2\pi i\cdot\sum_{m\equiv b-a(k),\atop m\equiv d-c(l)}
\frac{\exp(\frac{\pi i (k+l)}{kl}(\tau m^2+2mu)-2\pi i\frac{mw}{l}+
2\pi i (m\tau+u)(n_0-\frac{d}{l}+\frac{a}{k}))}{1-\exp(2\pi i(m\tau+u))}
\end{align*}
where $n_0$ is the minimal integer $n$ such that $n\ge\frac{w_2+d}{l}-
\frac{a}{k}$.  
It is easy to see that the function $F_{a,b,c}^d(u,w)$ is meromorphic in $u$
with poles of order 1 at the lattice points and 
satisfies (\ref{per1}) and (\ref{per2}).
Notice that the residues of $F_{a,b,c}^d(u,w)$ and
of $G_{a,b,c}^d(u,w)$ at $u=0$ are equal to 1. Therefore,
for fixed $w$, $a$, and $d$ the functions 
$H_{b,c}(u)=G_{a,b,c}^d(u,w)-F_{a,b,c}^d(u,w)$ are holomorphic
in $u$ and satisfy
\begin{equation}\label{per3}
H_{b,c}(u+1)=\exp(2\pi i(\frac{b}{k}-\frac{c}{l}))H_{b,c}(u)
\end{equation}
\begin{equation}\label{per4}
H_{b,c}(u+\tau)=\exp(-\frac{\pi i(k+l)}{kl}(2u+\tau)+2\pi i\frac{w}{l})
H_{b+1,c-1}(u)
\end{equation}
Let $r$ be the greatest common divisor of $k$ and $l$.
Iterating these equations $kl/r$ times we obtain
\begin{equation}\label{per5}
H_{b,c}(u+\frac{kl}{r})=H_{b,c}(u)
\end{equation}
\begin{equation}\label{per6}
H_{b,c}(u+\frac{kl}{r}\tau)=
\exp(-\pi i\frac{(k+l)kl}{r^2}\tau-2\pi i\frac{(k+l)u-kw}{r})
H_{b,c}(u)
\end{equation}
The space of holomorphic functions in $u$ satisfying (\ref{per5})
and (\ref{per6}) has the basis 
$e(s)=\th_{s/N}(\frac{(k+l)u-kw}{r},\frac{(k+l)kl}{r^2}\tau)$ where
$s\in \Z/N\Z$, $N=\frac{(k+l)kl}{r^2}$. 
Note that the action of the operator $f(u)\mapsto f(u+1)$ on
this basis is diagonal, hence, the space of holomorphic functions satisfying
(\ref{per3}) and (\ref{per6}) has as basis
$e((b/k-c/l+p)\frac{kl}{r})$, where
$p\in\Z/\frac{k+l}{n}\Z$. Thus, we should have
$$H_{b,c}(u)=\sum_p f_{b,c,p}e((b/k-c/l+p)\frac{kl}{r})$$ 
where $p$ runs through $\Z/\frac{k+l}{n}\Z$. 
Note that in this formula we have to fix representatives for 
the residue classes $b$ and $c$. A more convenient point of view is
the following. Consider the set $T$ of triples of integers
$(b,c,p)$ modulo the equivalence relation generated by 
\begin{align*}
&(b,c,p)\equiv (b,c, p+\frac{k+l}{r}),\\
&(b,c,p)\equiv (b+k,c,p-1),\\
&(b,c,p)\equiv (b,c+l,p+1)
\end{align*} 
Forgetting $p$ we get the projection 
$\phi_1:T\ra\Z/k\Z\times\Z/l\Z$. On the other hand, we have
a well-defined map
$$\phi_2:T\ra\Z/N\Z:(b,c,p)\mapsto (\frac{b}{k}-\frac{c}{l}+p)\frac{kl}{r}$$
Now the above formula can be rewritten as
$$H_{b,c}(u)=\sum_{\si\in\phi_1^{-1}(b,c)}f_{\si} e(\phi_2(\si)).$$
Now taking into account the relation (\ref{per4}) we derive the
relation $f_{b,c,p}=f_{b+1,c-1,p}$. Consider the action of $\Z$ 
on $T$ generated by $(b,c,p)\mapsto (b+1,c-1,p)$. It is easy
to see that the orbits of this action coincide with the fibers of
the map
$$\phi_3:T\ra\Z/(k+l)\Z:(b,c,p)\mapsto b+c+kp.$$
Therefore, the coefficients $f_{\si}$ depend only on $\phi_3(\si)$.
Summarizing (and recalling the dependence on $a$, $d$ and $w$) 
we can write
\begin{equation}\label{mainfun}
G_{a,b,c}^d(u,w)-F_{a,b,c}^d(u,w)=
\sum_{\si\in\phi_1^{-1}(b,c)}f_{a,\phi_3(\si)}^d(w) 
\th_{\phi_2(\si)/N}(\frac{(k+l)u-kw}{r},\frac{(k+l)kl}{r^2}\tau).
\end{equation}

\subsection{}\label{homot}
Recall that the addition formula for theta-functions
(see \cite{Theta}, prop. 6.4) gives the following identity
\begin{align*}
&\th_{b/k}(kx+u,k\tau)\th_{c/l}(lx+v,l\tau)=\\
&\sum_{p\in\Z/\frac{k+l}{r}\Z}
\th_{(\frac{b}{k}-\frac{c}{l}+p)\frac{r}{k+l}}(
\frac{(k+l)u-kw}{r},\frac{(k+l)kl}{r^2}\tau)
\th_{\frac{b+c+kp}{k+l}}((k+l)x+w,(k+l)\tau)
\end{align*}
where $w=u+v$. In this formula one has to fix representatives of
the residue classes $b$ and $c$. Using the notation above we
can rewrite it as follows
\begin{align*}
&\th_{b/k}(kx+u,k\tau)\th_{c/l}(lx+v,l\tau)=\\
&\sum_{\si\in\phi_1^{-1}(b,c)}
\th_{\frac{\phi_2(\si)}{N}}(\frac{(k+l)u-kw}{r},\frac{(k+l)kl}{r^2}\tau)
\th_{\frac{\phi_3(\si)}{k+l}}((k+l)x+w,(k+l)\tau)
\end{align*}
Comparing this with (\ref{mainfun}) we conclude that we can write
the difference between the two types of triple product of
the elements 
$\a\in H^1(L(-k,0))$, $\b_1\in H^0(L(k,u))$ and 
$\b_2\in H^0(L(l,v))$ in the form
$$m_3(\a,\b_1,\b_2)-m_3'(\a,\b_1,\b_2)=
n_2(\a,\b_1\b_2)$$
where the linear map $n_2:H^1(L(-k,0))\otimes H^0(L(k+l,w))\ra
H^0(L(l,w))$ is defined by the formula
\begin{equation}
n_2(\a,\th_{\frac{q}{k+l}}((k+l)x+w,(k+l)\tau))=
\sum_{d\in\Z/l\Z}f_{a,q}^d \th_{d/l}(lx+w,l\tau).
\end{equation}
where $q\in\Z/(k+l)\Z$.

Note that we have
$$m_3(\b_2,\b_1,\a)=-\pr(\b_2(\dbar^{-1}(\b_1\a))=
-\pr(\dbar^{-1}(\a\b_1)\b_2)=-m_3(\a,\b_1,\b_2).$$
Similar formula holds for $m'_3$ so we get a homotopy
between $m_3(\b_2,\b_1,\a)$ and $m'_3(\b_2,\b_1,\a)$
by setting $n_2(\b,\a)=n_2(\a,\b)$.

\subsection{}
We claim that using $A_{\infty}$-axioms one can
express all transversal triple products between line bundles on $E_{\tau}$
in terms of the ones we computed in \ref{comp} and in terms of the
triple products of the form
$$H^0(L(d,u_1))\otimes H^1(L(-d,u_2))\otimes H^0(L(d,u_3))\ra
H^0(L(d,u_1+u_2+u_3))$$
where $u_1+u_2\not\in\Z+\Z\tau$, $u_2+u_3\not\in\Z+\Z\tau$.
The latter triple products in both categories coincide as the univalued
well-defined Massey products (see \cite{P}).

First consider the triple product $m_3(r,s,t)$ where
$r\in H^1(L(-d_1,u_1)$, $s\in H^0(L(d_2,u_2))$, $t\in H^0(L(d_3,u_3))$,
where $d_i>0$ and $d_2+d_3>d_1$. If $d_2=d_1$ and $u_1+u_2\not\in\Z+\Z\tau$
this product
reduces to the one considered in \ref{comp} by a translation of
argument. Otherwise, there are two possibilities:

\noindent a) $d_2>d_1$. In this case we can write $s$ as a linear
combination of products $s's''$ where
$s'\in H^0(L(d_1,u'_2))$, $s''\in H^0(L(d_2-d_1,u''_2))$,
$u'_2+u''_2=u_2$. Moreover, we can assume that $u_1+u'_2$ doesn't
belong to the lattice $\Z+\Z\tau$. Then applying the 
$A_{\infty}$-constraint to the quadruple $r,s',s'',t$ we obtain that
$m_3(r,s's'',t)$ is a linear combination of $m_3(r,s',s'')t$ and
$m_3(r,s',s''t)$ (note that $rs'=0$ by assumption while
$m_3(s',s'',t)=0$ as an element of $H^{-1}(L(d_2+d_3,u_2+u_3))$).

\noindent b) $d_2<d_1$. In this case we write $t$ as a linear
combination of products $t't''$ where
$t'\in H^0(L(d_1-d_2,u'_3)$, $t''\in H^0(L(d_2+d_3-d_1, u''_3))$,
$u_3=u'_3+u''_3$. Moreover, we can assume that $u_1+u_2+u'_3$
doesn't belong to the lattice $\Z+\Z\tau$. This implies
that $m_3(r,s,t')=0$. Hence, applying the $A_{\infty}$-constraint
to the quadruple $r,s,t',t''$ we get
$m_3(r,s,t't'')=\pm m_3(rs,t',t'')\pm m_3(r,st',t'')$ (note that 
$m_3(s,t',t'')=0$).

One deals similarly with transversal products of the type
$$H^0(L(d_1,u_1))\otimes H^0(L(d_2,u_2))\otimes H^1(L(-d_3,u_3))\ra
H^0(L(d_1+d_2-d_3,u_1+u_2+u_3))$$
where $d_i>0$.

Now let us consider $m_3(r,s,t)$ where $r\in H^0(L(d_1,u_1))$,
$s\in H^1(L(-d_2,u_2))$, $t\in H^0(L(d_3,u_3))$,
$d_1+d_3>d_2$, $d_i>0$. If $d_1=d_3=1$ then we necessarily have
$d_2=1$, hence this is a univalued Massey product. Otherwise,
we can use recursion in $max(d_1,d_3)$. Indeed, if $d_1>1$ then
we can write $r$ as a linear combination of $r'r''$ where
$r'\in H^0(L(1,u'_1))$, $r''\in H^0(L(d_1-1,u''_2))$.
Now
$$m_3(r'r'',s,t)=\pm r'm_3(r'',s,t) \pm m_3(r',r''s,t) \pm
m_3(r',r'',st) \pm m_3(r',r'',s)t$$
where the products $m_3(r',r'',st)$ and $m_3(r',r'',s)$ are of the type
treated above. Similarly, if $d_3>1$ then we can write $t$ as
linear combination of $t't''$ to get a recursion formula for
$m_3(r,s,t)$.

Finally, using the compatibility of $m_3$ with the Serre duality
and the obvious cyclic symmetry of $m'_3$ we can reduce all non-zero
transversal higher products $m_3$ and $m'_3$ to the ones considered above.
For example, the product
$$H^1(L_1)\otimes H^1(L_2)\otimes H^0(L_3)\ra H^1(L_1\otimes L_2\otimes L_3)
$$
by Serre duality reduces to
$$H^1(L_2)\otimes H^0(L_3)\otimes H^0((L_1\otimes L_2\otimes L_3)^{-1})\ra
H^0(L_1^{-1}).$$

Using these formulas it is easy to see that the homotopy constructed
in \ref{homot} induces homotopy in all the other cases of
transversal triple products between the line bundles.

%

\vspace{3mm}

{\sc Department of Mathematics, Harvard University, Cambridge,
 MA 02138}

{\it E-mail address:} apolish@@math.harvard.edu

\end{document}